# Pantazi's Theorem Regarding the Bi-orthological Triangles

Prof. Ion Pătraşcu, The National College "Fraţii Buzeşti", Craiova, Romania
Prof. Florentin Smarandache, University of New Mexico, U.S.A.

In this article we'll present an elementary proof of a theorem of Alexandru Pantazi (1896-1948), Romanian mathematician, regarding the bi-orthological triangles.

1. **Orthological triangles**

**Definition**
The triangle $ABC$ is orthologic in rapport to the triangle $A_1B_1C_1$ if the perpendiculars constructed from $A, B, C$ respectively on $B_1C_1, C_1A_1$ and $A_1B_1$ are concurrent. The concurrency point is called the orthology center of the triangle $ABC$ in rapport to triangle $A_1B_1C_1$.

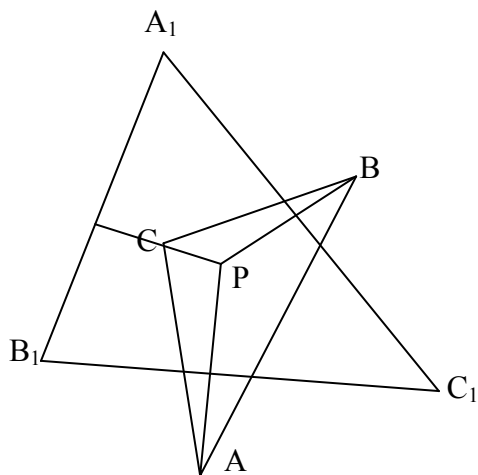

Fig. 1

In figure 1 the triangle $ABC$ is orthologic in rapport with $A_1B_1C_1$, and the orthology center is $P$.

2. **Examples**

a) The triangle $ABC$ and its complementary triangle $A_1B_1C_1$ (formed by the sides' middle) are orthological, the orthology center being the orthocenter $H$ of the triangle $ABC$.
Indeed, because $B_1C_1$ is a middle line in the triangle $ABC$, the perpendicular from $A$ on $B_1C_1$ will be the height from $A$. Similarly the perpendicular from $B$ on $C_1A_1$ and the perpendicular from $C$ on $A_1B_1$ are heights in $ABC$, therefore concurrent in $H$ (see Fig. 2)

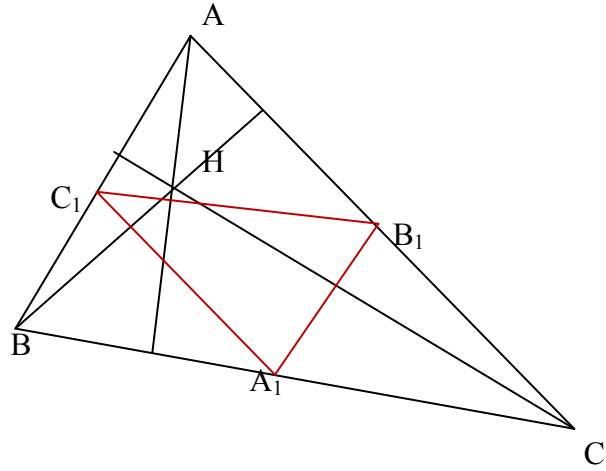

Fig. 2

b) **Definition**

Let $D$ a point in the plane of triangle $ABC$. We call the circum-pedal triangle (or meta-harmonic) of the point $D$ in rapport to the triangle $ABC$, the triangle $A_1B_1C_1$ of whose vertexes are intersection points of the Cevianes $AD, BD, CD$ with the circumscribed circle of the triangle $ABC$.

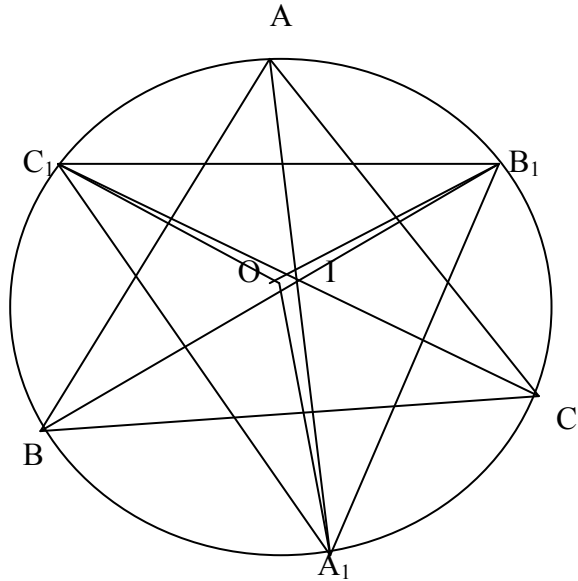

Fig. 3

The triangle circum-pedal $A_1B_1C_1$ of the center of the inscribed circle in the triangle $ABC$ and the triangle $ABC$ are orthological (Fig. 3).

The points $A_1, B_1, C_1$ are the midpoints of the arcs $\widehat{BC}, \widehat{CA}$ respectively $\widehat{AB}$. We have $\widehat{A_1B} \equiv \widehat{A_1C}$, it results that $A_1B = A_1C$, therefore $A_1$ is on the perpendicular bisector of $BC$, and therefore the perpendicular raised from $A_1$ on $BC$ passes through $O$, the center of the circumscribed circle to triangle $ABC$. Similarly the perpendiculars raised from $B_1, C_1$ on $AC$ respectively $AB$ pass through $O$. The orthology center of triangle $A_1B_1C_1$ in rapport to $ABC$ is $O$

### 3. The characteristics of the orthology property

The following Lemma gives us a necessary and sufficient condition for the triangle $ABC$ to be orthologic in rapport to the triangle $A_1B_1C_1$.

**Lemma**

The triangle $ABC$ is orthologic in rapport with the triangle $A_1B_1C_1$ if and only if:

$$\overrightarrow{MA} \cdot \overrightarrow{B_1C_1} + \overrightarrow{MB} \cdot \overrightarrow{C_1A_1} + \overrightarrow{MC} \cdot \overrightarrow{A_1B_1} = 0 \tag{1}$$

for any point $M$ from plane.

**Proof**

In a first stage we prove that the relation from the left side, which we'll note $E(M)$ is independent of the point $M$.

Let $N \neq M$ and $E(N) = \overrightarrow{NA} \cdot \overrightarrow{B_1C_1} + \overrightarrow{NB} \cdot \overrightarrow{C_1A_1} + \overrightarrow{NC} \cdot \overrightarrow{A_1B_1}$

Compute $E(M) - E(N) = \overrightarrow{MN} \cdot \left(\overrightarrow{BC} + \overrightarrow{CA} + \overrightarrow{AB}\right)$.

Because $\overrightarrow{BC} + \overrightarrow{CA} + \overrightarrow{AB} = 0$ we have that $E(M) - E(N) = \overrightarrow{MN} \cdot \vec{0} = 0$.

If the triangle $ABC$ is orthologic in rapport to $A_1B_1C_1$, we consider $M$ their orthologic center, it is obvious that (1) is verified. If (1) is verified for a one point, we proved that it is verified for any other point from plane.

Reciprocally, if (1) is verified for any point $M$, we consider the point $M$ as being the intersection of the perpendicular constructed from $A$ on $B_1C_1$ with the perpendicular constructed from $B$ on $C_1A_1$. Then (1) is reduced to $\overrightarrow{MC} \cdot \overrightarrow{A_1B_1} = 0$, which shows that the perpendicular constructed from $C$ on $\overrightarrow{A_1B_1}$ passes through $M$. Consequently, the triangle $ABC$ is orthologic in rapport to the triangle $A_1B_1C_1$.

### 4. The symmetry of the orthology relation of triangles

It is natural to question ourselves that given the triangles $ABC$ and $A_1B_1C_1$ such that $ABC$ is orthologic in rapport to $A_1B_1C_1$, what are the conditions in which the triangle $A_1B_1C_1$ is orthologic in rapport to the triangle $ABC$.

The answer is given by the following

**Theorem** (The relation of orthology of triangles is symmetric)

If the triangle $ABC$ is othologic in rapport with the triangle $A_1B_1C_1$ then the triangle $A_1B_1C_1$ is also orthologic in rapport with the triangle $ABC$.

**Proof**

We'll use the lemma. If the triangle $ABC$ is orthologic in rapport with $A_1B_1C_1$ then

$$\overrightarrow{MA} \cdot \overrightarrow{B_1C_1} + \overrightarrow{MB} \cdot \overrightarrow{C_1A_1} + \overrightarrow{MC} \cdot \overrightarrow{A_1B_1} = 0$$

for any point $M$. We consider $M = A$, then we have

$$\overrightarrow{AA} \cdot \overrightarrow{B_1C_1} + \overrightarrow{AB} \cdot \overrightarrow{C_1A_1} + \overrightarrow{AC} \cdot \overrightarrow{A_1B_1} = 0.$$

This expression is equivalent with

$$\overrightarrow{A_1A} \cdot \overrightarrow{BC} + \overrightarrow{A_1B_1} \cdot \overrightarrow{CA} + \overrightarrow{A_1C_1} \cdot \overrightarrow{AB} = 0$$

That is with (1) in which $M = A_1$, which shows that the triangle $A_1B_1C_1$ is orthologic in rapport to triangle $ABC$.

### Remarks
1. We say that the triangles $ABC$ and $A_1B_1C_1$ are orthological if one of the triangle is orthologic in rapport to the other.
2. The orthology centers of two triangles are, in general, distinct points.
3. The second orthology center of the triangles from a) is the center of the circumscribed circle of triangle $ABC$.
4. The orthology relation of triangles is reflexive. Indeed, if we consider a triangle, we can say that it is orthologic in rapport with itself because the perpendiculars constructed from $A, B, C$ respectively on $BC, CA, AB$ are its heights and these are concurrent in the orthocenter $H$.

### 5. Bi-orthologic triangles

### Definition
If the triangle $ABC$ is simultaneously orthologic to triangle $A_1B_1C_1$ and to triangle $B_1C_1A_1$, we say that the triangles $ABC$ and $A_1B_1C_1$ are bi-orthologic.

### Pantazi's Theorem
If a triangle $ABC$ is simultaneously orthologic to triangle $A_1B_1C_1$ and $B_1C_1A_1$, then the triangle $ABC$ is orthologic also with the triangle $C_1A_1B_1$.

### Proof
Let triangle $ABC$ simultaneously orthologic to $A_1B_1C_1$ and to $B_1C_1A_1$, using lemma, it results that

$$\overrightarrow{MA} \cdot \overrightarrow{B_1C_1} + \overrightarrow{MB} \cdot \overrightarrow{C_1A_1} + \overrightarrow{MC} \cdot \overrightarrow{A_1B_1} = 0 \tag{2}$$
$$\overrightarrow{MA} \cdot \overrightarrow{C_1A_1} + \overrightarrow{MB} \cdot \overrightarrow{A_1B_1} + \overrightarrow{MC} \cdot \overrightarrow{B_1C_1} = 0 \tag{3}$$

For any $M$ from plane.

Adding the relations (2) and (3) side by side, we have:

$$\overrightarrow{MA} \cdot \left(\overrightarrow{B_1C_1} + \overrightarrow{C_1A_1}\right) + \overrightarrow{MB} \cdot \left(\overrightarrow{C_1A_1} + \overrightarrow{A_1B_1}\right) + \overrightarrow{MC} \cdot \left(\overrightarrow{A_1B_1} + \overrightarrow{B_1C_1}\right) = 0$$

Because

$$\overrightarrow{B_1C_1} + \overrightarrow{C_1A_1} = \overrightarrow{B_1A_1},\ \overrightarrow{C_1A_1} + \overrightarrow{A_1B_1} = \overrightarrow{C_1B_1},\ \overrightarrow{A_1B_1} + \overrightarrow{B_1C_1} = \overrightarrow{A_1C_1}$$

(Chasles relation), we have:

$$\overrightarrow{MA} \cdot \overrightarrow{B_1A_1} + \overrightarrow{MB} \cdot \overrightarrow{C_1B_1} + \overrightarrow{MC} \cdot \overrightarrow{A_1C_1} = 0$$

for any $M$ from plane, which shows that the triangle $ABC$ is orthologic with the triangle $C_1A_1B_1$ and the Pantazi's theorem is proved.

**Remark**

The Pantazi's theorem can be formulated also as follows: If two triangles are bi-orthologic then these are tri-orthologic.

**Open Questions**
1) Is it possible to extend Pantazi's Theorem (in 2D-space) in the sense that if two triangles $A_1B_1C_1$ and $A_2B_2C_2$ are bi-orthological, then they are also *k*-orthological, where *k = 4, 5,* or *6?*
2) Is it true a similar theorem as Pantazi's for two bi-homological triangles and bi-orthohomological triangles (in 2D-space)? We mean, if two triangles $A_1B_1C_1$ and $A_2B_2C_2$ are bi-homological (respectively bi-orthohomological), then they are also *k*-homological (respectively *k*-orthohomological), where *k = 4, 5,* or *6?*
3) How the Pantazi Theorem behaves if the two bi-orthological non-coplanar triangles $A_1B_1C_1$ and $A_2B_2C_2$ (if any) are in the 3D-space?
4) Is it true a similar theorem as Pantazi's for two bi-homological (respectively bi-orthohomological) non-coplanar triangles $A_1B_1C_1$ and $A_2B_2C_2$ (if any) in the 3D-space?
5) Similar questions as above for bi-orthological / bi-homological / bi-orthohomological polygons (if any) in 2D-space, and respectively in 3D-space.
6) Similar questions for bi-orthological / bi-homological / bi-orthohomological polyhedrons (if any) in 3D-space.

**References**

[1]    Cătălin Barbu – Teoreme fundamentale din geometria triunghiului, Editura Unique, Bacău, 2007.
[2]    http://garciacapitan.auna.com/ortologicos/ortologicos.pdf.5